\documentclass[12pt]{amsart}
\usepackage{graphics}

\newcommand {\supplus}{\mathop{{\supset}\llap{\raise 
0.5pt\hbox{\normalfont\small+}\hskip 0.5pt}}} 

\newcommand {\subplus}{\mathop{{\subset}\llap{\raise 
0.5pt\hbox{\normalfont\small+}\hskip 0.5pt}}}  

\newcommand {\Cee}    {{\mathbb  C}}

\newcommand {\Zee}    {{\mathbb  Z}}

\newcommand {\fas}    {{\mathfrak{as}}}
 
\newcommand {\fb}     {{\mathfrak{b}}}

\newcommand {\fd}     {{\mathfrak{d}}}

\newcommand {\fg}     {{\mathfrak{g}}}    %
\newcommand {\fgl}    {{\mathfrak{gl}}}  %
\newcommand {\fh}     {{\mathfrak{h}}}

\newcommand {\fk}     {{\mathfrak{k}}}

\newcommand {\fle}    {{\mathfrak{le}}}
\newcommand {\fm}     {{\mathfrak{m}}}
\newcommand {\fN}     {{\mathfrak{N}}}  %
\newcommand {\fn}     {{\mathfrak{n}}}

\newcommand {\fo}     {{\mathfrak{o}}}
\newcommand {\fosp}   {{\mathfrak{osp}}}
\newcommand {\fp}    {{\mathfrak{p}}}   %
\newcommand {\fpe}    {{\mathfrak{pe}}}   %

\newcommand {\fpo}    {{\mathfrak{po}}}

\newcommand {\fpsl}   {{\mathfrak{psl}}}

\newcommand {\fpsq}   {{\mathfrak{psq}}}
\newcommand {\fq}     {{\mathfrak{q}}}

\newcommand {\fs}     {{\mathfrak{s}}}

\newcommand {\fsh}    {{\mathfrak{sh}}}

\newcommand {\fsl}    {{\mathfrak{sl}}}
\newcommand {\fsle}   {{\mathfrak{sle}}}
\newcommand {\fsm}    {{\mathfrak{sm}}}
\newcommand {\fsp}    {{\mathfrak{sp}}}
\newcommand {\fspe}   {{\mathfrak{spe}}}

\newcommand {\fsvect} {{\mathfrak{svect}}}

\newcommand {\fvect}  {{\mathfrak{vect}}}   %

\newcommand {\cal} {\mathcal}

\def \opname#1#2%
  {\expandafter\newcommand \csname #1\endcsname {{\mathop{#2}\nolimits}}}

\newcommand{\rmname}[1]
  {\expandafter\newcommand \csname #1\endcsname {{\operatorname{#1}}}}

\newcommand{\rmnameii}[2]
  {\expandafter\newcommand \csname #1\endcsname {{\operatorname{#2}}}}

\rmname{act}
\rmname{Ad}
\rmname{Add}
\rmname{ad}
\rmname{Alt}
\rmname{alt}
\rmname{Ann}
\rmname{antidiag}
\rmname{Ber}
\rmname{ber}
\rmname{Br}
\rmname{card}
\rmname{ch}
\rmname{Char}
\rmname{cem}
\rmname{cj}
\rmname{Cliff}
\rmname{cntr}
\rmname{codim}
\rmname{coind}
\rmname{const}
\rmname{col}
\rmname{cork}
\rmname{cpr}
\rmname{diag}
\opname{Div}{div}
\rmname{Def}
\rmname{Der}
\rmname{Dim}
\rmname{End}
\rmname{Even}
\rmname{Ext}
\rmname{gr}
\rmname{Hom}
\rmname{HT}
\rmnameii{Ht}{ht}
\rmname{hwt}
\rmname{Id}
\rmname{id}
\rmname{ind}
\rmname{Ind}
\rmname{Inf}
\rmname{irr}
\rmname{Le}
\rmname{Lie}
\rmname{lwt}
\rmname{mult}
\rmname{Mor}
\rmname{nm}
\rmname{Ob}
\rmname{Odd}
\rmname{Osc}
\rmname{per}
\rmname{Pic}
\rmname{pr}
\rmname{pro}
\rmname{Prime}
\rmname{Proj}
\rmname{prt}
\rmname{pt}
\rmname{Q}
\rmname{qet}
\rmname{qtr}
\rmname{rd}
\rmname{rk}
\rmname{row}
\rmname{Res}
\rmname{salt}
\rmname{Sch}
\rmname{SBr}
\rmname{scalar}
\rmname{Ser}
\rmname{sign}
\rmname{Smbl}
\rmname{spin}
\rmname{ssym}
\rmname{str}
\rmname{st}
\rmname{sgn}
\rmname{sq}
\rmname{symm}
\rmname{supp}
\rmname{Supp}
\rmname{St}
\rmname{Spec}
\rmname{Spm}
\rmname{tr}
\rmname{vpt}
\rmname{weyl}
\rmname{Weyl}
\rmname{Witt}

\opname{vvol}  {{v\hspace{-0.1ex}o\hspace{-0.02ex}l\/}}
\opname{pnt}  {\text{\normalfont pt}}
\opname{Span} {{Span}}
\opname{slim} {\overline{\lim}}
\opname{Vol}  {{V\hspace{-0.55ex}o\hspace{-0.02ex}l\/}}
\opname{Par}  {{P\hspace{-0.3ex}a\hspace{-0.05ex}r\/}}

\rmname{Mat}
\rmname{Bil}
\rmname{Diff}
\rmname{Ker}
\rmname{Herm}
\rmname{Coker}
\rmname{Conn}
\rmname{Covect}
\rmname{Vect}
\rmname{Int}

\rmnameii {IM} {Im}
\rmnameii {RE} {Re}

\opname{Aut} {{A\hspace{-0.2ex}u\hspace{-0.1ex}t\/}}
\opname{GL} {{G\hspace{-0.3ex}L}}
\opname{SL} {{S\hspace{-0.3ex}L}}
\opname{Exp} {{E\hspace{-0.2ex}x\hspace{-0.1ex}p\/}}
\opname{GQ} {{G\hspace{-0.2ex}Q}}
\opname{OSp} {{O\hspace{-0.25ex}S\hspace{-0.15ex}p\/}}
\opname{Out} {{O\hspace{-0.25ex}u\hspace{-0.15ex}t\/}}
\opname{Spp} {{S\hspace{-0.2ex}p\/}}
\opname{SpO} {{S\hspace{-0.2ex}p\hspace{-0.02ex}O\/}}
\opname{Pe} {{P\hspace{-0.25ex}e\/}}
\opname{SPe} {{S\hspace{-0.25ex}P\hspace{-0.25ex}e\/}}
\opname{Spin} {{S\hspace{-0.25ex}p\hspace{-0.05ex}i\hspace{-0.1ex}n\/}}
\opname{Iso} {{I\hspace{-0.25ex}s\hspace{-0.1ex}o\/}}
\opname{SSPe} {{S\hspace{-0.25ex}S\hspace{-0.15ex}P\hspace{-0.25ex}e\/}}
\opname{PeU} {{P\hspace{-0.25ex}e\hspace{-0.1ex}U\/}}
\opname{QU} {{Q\hspace{-0.15ex}U\/}}
\opname{U} {{U\/}}

\opname{cGQ} {{\cal G \hspace{-0.2em} Q \/}}
\opname{cSL} {{\cal S \hspace{-0.2em} L \/}}
\opname{cGL} {{\cal G \hspace{-0.2em} L \/}}
\opname{cOSp} {{\cal O \hspace{-0.2em} S \hspace{-0.3em} \it p\/}}
\opname{cPe} {{\cal P \hspace{-1.5pt} \it e\/}}
\opname{cVect} {{\cal V \hspace{-1.5pt} \it e\hspace{-0.1ex}c\hspace{-0.1ex}t\/}}
\opname{cVol} {{\cal V \hspace{-1.5pt} \it o\hspace{-0.1ex}l\/}}
\opname{cAut} {{\cal A \hspace{-0.2em} \it u\hspace{-0.1em}t\/}}
\opname{cCovect} {{\cal C \hspace{-1.5pt}
     \it o\hspace{-0.1ex}v\hspace{-0.1ex}e\hspace{-0.1ex}c\hspace{-0.1ex}t\/}}
\opname{CW} {{C\hspace{-0.15ex}W}}

\newcommand {\eps} {\varepsilon}

\newcommand {\bcdot}   {\mathbin{\hbox{\raise.4ex\hbox{\bf.}}}} 

\newcommand {\secno} {}

\newtheorem{Theorem}{\secno Theorem}

\newtheorem{Lemma}[Theorem]{\secno Lemma}

\newenvironment {th*}[1]
    {\gdef\thname{#1} \begin{thn}}%
    {\end{thn}}
\newtheorem{thn}[Theorem] {\thname}

\theoremstyle{definition}

\newenvironment {ex*}[1]
    {\gdef\thname{#1} \begin{exn}}%
    {\end{exn}}
\newtheorem{exn}[Theorem]{\thname}

\theoremstyle{remark}

\newenvironment {rem*}[1]
    {\gdef\thname{#1} \begin{remn}}%
    {\end{remn}}
\newtheorem{remn}[Theorem]{\thname}

\newcommand {\ssec}{\subsection*}

\begin{document}
	
\title[Defining 
relations for classical Lie superalgebras]{Defining 
relations for classical Lie superalgebras without Cartan matrices}

\author{P.~Grozman${}^1$, D.~Leites${}^1$; 
E.~Poletaeva${}^2$}

\address{${}^1$(Correspondence): Department of
Mathematics, University of Stockholm, Kr\"aftriket hus
6, SE-106 91, Stockholm, Sweden; e-mail: mleites@ matematik.su.se;
${}^2$Department of Mathematics, Lund University, Sweden;
elena@maths.lth.se}

\thanks{We were financially supported: P.G. by the Swedish Institute 
in 1993--95; D.L. by I.~Bendixson grant, an NSF grant via IAS, 
Princeton, in 1989 and by NFR during 1986--1999; D.L. and E.P. by 
SFB-170 in June--July of 1990. We are particularly thankful to 
Yu.~Kochetkov, actually a co-author, for calculating the relations 
for $\fsh(0|n)$.}

\begin{abstract}
The analogs of Chevalley generators are offered for simple (and close
to them) $\Zee$-graded complex Lie algebras and Lie superalgebras of
polynomial growth without Cartan matrix.  We show how to derive the
defining relations between these generators and explicitly write them
for a \lq\lq most natural" (\lq\lq distinguished" in terms of Penkov
and Serganova) system of simple roots.  The results are given mainly
for Lie superalgebras whose component of degree zero is a Lie algebra
(other cases being left to the reader).  Observe
presentations of exceptional Lie superalgebras and Lie superalgebras 
of hamiltonian vector fields. Now we can, at last, $q$-quantize the Lie
Lie superalgebras of hamiltonian vector fields and Poisson superalgebras.
\end{abstract}

\keywords{Lie superalgebras, defining relations.}

\subjclass{17A70 (Primary) 17B01, 17B70 (Secondary)}

\maketitle

\section*{\protect \S 1. Preliminaries}

After Berezin formulated in 1971 to one of us (DL) the problem which
in modern terms would have sounded ``define supermanifolds and
construct differential supergeometry'' a first step was to look after
examples of Lie superalgebras that naturally appear in mathematics. 
Homotopy rings with respect to Whitehead's product are examples of
such Lie rings, but nobody, so far, described these natural rings in
any case.  By description we mean identification of the semisimple
part and radical.  A reason for such careless attention to these most
natural rings becomes clear soon after one tackles the problem: they
are nilpotent, hence, not so interesting in a sense (simple Lie
(super)algebras have a richer structure and interesting representation
theory).  A result of C.~L\"ofvall and J.-E.~Roos indicates that this
deduction might be too hasty.

In a similar problem C.~L\"ofvall and J.-E.~Roos \cite{LR} made a
break-through observation: they not only found traces left by simple
Lie superalgebras where nothing indicated them but also identified
these superalgebras as a ``positive part'' of certain twisted loops
algebras with values in simple Lie superalgebras.  The paper \cite{LR}
is, clearly, the first in a series to appear, where, among other
things, L\"ofvall and Roos will need presentations (in other words,
generators and defining relations) of the {\it positive} part
$\fg_{>}=\mathop{\sum}\limits_{i>0}\fg_{i}$ of certain (twisted) loop
superalgebra $\fg=\mathop{\sum}\limits_{i>0}\fg_{i}$, associated with
the simple (or a close to the simple) finite dimensional Lie
superalgebras.

The results obtained here (an extension of \cite{GLP}), and those of
\cite{GL1}, \cite{GL2}, \cite{Y}, \cite{LSe}, as well as those which,
though listed as open problems, are supplied with an instruction how
to solve them, describe how to shorthand the presentation needed, both
generators and relations, compare with appalling presentation of
\cite{T} or implicit presentations of the positive parts of vectorial
algebras in \cite{FF}.

We hope that our results contribute to the fest on the occasion of
Jan-Erik's birthday and make his calculations, if not life, easier.

Namely, we observe that even in the absence of $\fg_{0}$ there is a
concise way to encode the presentation.  Indeed, in the majority of
cases $\fg_{>}$ is the direct sum of irreducible $\fg_{0}$-modules
$\fg_{i}$, and, as Lie superalgebra, $\fg_{>}$ is generated by
$\fg_{1}$.  If $\fg_{1}$ is irreducible, as in cases of L\"ofvall and
Roos, then, instead of $\dim~\fg_{1}$ generators, it suffices to take just
one, any one vector, say the lowest weight vector. 

The space of relations (same as the space of generators in the general
case) must not split into the direct sum of {\it irreducible}
$\fg_{0}$-modules, but, nevertheless, one can list only the vacuum
vectors, i.e., the lowest AND highest weight vectors (since some
modules can be glued in indecomposable conglomerates, we need both). 
{\it Mathematica}-based package SuperLie (\cite{G}) helps to find these
vacuum vectors.

To reduce volume of the paper, we did not reproduce the standard
homological interpretation of relations and spectral sequence leading
to the answer; for Lie algebras it is expressed in \cite{LP} and its
superization is straightforward via the Sign Rule.  Observe only that
since relations represent homology class, they can be ``pure'' or
``dirty'' if defined modulo boundaries.

{\bf History and an overview}.  The traditional way to determine
classical simple finite dimensional Lie algebra (over $\Cee$) is via
Chevalley generators, though other generators are possible.  For
discussion of other possibilities with examples see \cite{GL1}.  Recently a
presentation of simple Lie superalgebras of the four Cartan series of
vector fields was given \cite{LP} and, together with Serre relations
for affine Kac--Moody algebras \cite{K1}, this completed description
of presentations of simple $\Zee$-graded Lie algebras of polynomial
growth (presentations with respect to other choices of generators are certainly 
possible).

Here we consider simple $\Zee$-graded Lie
superalgebras of polynomial growth (and close to them \lq\lq
classical" Lie superalgebras, such as central extensions of the simple
ones, their algebras of differentiations, etc.).  Their list is
conjecturally (\cite{LS1}) completed and consists of 

$\bullet$ the finite dimensional ones (classification results by Kaplansky
and Nahm-Rittenberg-Scheunert \cite{FK}, \cite{NRS} were skillfully
rounded up by Kac \cite{K2}),

$\bullet$ the vectorial algebras, i.e., algebras of vector fields,
(classification announced \cite{LS1} and partly proved \cite{LS2} by
Leites and Schepochkina; the proof was again quickly rounded up by
Kac and Cheng \cite{K3}, bar some gaps, see \cite{Sh5}),

$\bullet$ the twisted loop algebras (with symmetrizable Cartan matrix
\cite{vdL}), or obtained as twisted loops \cite{FLS}, and

$\bullet$ the stringy (i.e., vectorial algebras pertaining to string
theories) algebras (for their intrinsic definition and list see
\cite{GLS1}).

In terms of presentations, another subdivision is more natural: 

(a) the algebras of the form $\fg(A)$ with Cartan matrix $A$ (subdivided
into subclasses (a${}_{s}$) with symmetrizable Cartan matrix and
(a${}_{n}$) with non-symmetrizable Cartan matrix, 

(q) the series $\fp\fs\fq$ and its relatives (central extensions,
exterior differentiations, etc.)  and

(v) the vectorial algebras and their relatives, the members of the
subclass are easily recognized by lack of the property ``if $\alpha$
is a root, then so is $-\alpha$''.

For $\fg(A)$ (both subcases) a very redundant presentation is given in
\cite{LSS} and, the minimal one, in \cite{GL2}; their $q$-quantization
(for symmetrizable $A$) is described in \cite{Y}.  The redundant
presentation \cite{LSS}, though very long, has an advantage: it only
involves Serre relations.  Regrettably, it is so redundant that
practically it is useless.

For $\fq(n)$ and $\fq(n)^{(1)}$ see \cite{LSe}; presentation of 
twisted loops is an open problem.

For presentations of vectorial Lie algebras see \cite{LP}; the series
with $\fg_{0}$ a Lie algebra see \cite{GLP}; here we also consider
several cases left in \cite{GLP} as open problems and the Lie
superalgebras of Hamiltonian vector fields and its central extension:
Poisson superalgebra.  The last two cases are of interest in relation
with spinor-oscillator representations and its $q$-quantization, see 
\cite{Kl}, \cite{LSh}. 

Though several results describing presentations of simple and close to
them vectorial superalgebras were obtained a while ago (\cite{Ko},
\cite{T}) they are given in the form too bulky to grasp or implicit
(\cite{U}; \cite{FF}).  A simplification of presentations is
desirable: for $\fg_{0}=\fsh(0|n)$ the dimension of the space of
relations computed in \cite{T} for $n=5$ is equal to 420 and grows
with $n$, whereas the total number of vacum vectors in the space of
relations is $<10$ and does not grow with $n$, cf.  Tables 2.1.2.

{\bf Problem formulation}.  Consider $\Zee$-graded Lie superalgebras
$\fg= \mathop{\oplus}\limits_{i\in \Zee}\fg_i$ of the following two
types:

(1) vectorial Lie superalgebras, i.e., of finite {\it depth} $d$ (in
the above sum $i\geq d$), cf.  \cite{LP};

(2) of infinite depth, but not of the form $\fg(A)$ (the algebras 
$\fg(A)$ being already considered in complementing each other papers 
\cite{GL2} and \cite{Y}) or of type $\fq$ (for whose presentation see 
\cite{LSe}).

Among these algebras we will first consider the ones for which $\fg_0$
is a Lie algebra.  The general case is an open problem, which can be
solved any time for any given $\fg$ via the lines indicated here and
with the help of Grozman's SuperLie package.  

Observe that in \cite{GL2} and \cite{Y} all bases (systems of simple
roots or, rather, corresponding generators) are considered.  For the
vectorial algebras and superalgebras and for loop algebras with values
in vectorial superalgebras we have considered below just one of the possible
bases.  It can well happen that presentations corresponding to some
other base is nicer in some sense: e.g., for a rank $n$ simple Lie
algebra, Serre relations corresponding to $3n$ Chevalley generators
though numerous ($\sim n^2$) are very simple and easy to compute, unlike
a handful of independent on $n$ but more intricate relations between the pair of
``Jacobson generators'' considered in \cite{GL1}. Nevertheless, both 
presentations are needed. For a method of passage from base to base 
(an analog of the Weyl group) see \cite{PS}.

Let $\fg_{+} =\mathop{\oplus}\limits_{i>0}\fg_i$ and $\fg_{-} 
=\mathop{\oplus}\limits_{i<0}\fg_i$; let $\fn_{\pm}$ be a maximal 
nilpotent subalgebra of $\fg_0$ described in textbooks (e.g., 
\cite{OV}) if $\fg_0$ is a Lie algebra or in \cite{PS} (see also refs.  
therein) if $\fg_0$ is a Lie superalgebra.  We decompose $\fg$ into 
the sum $\fN_{-}\oplus \fh\oplus \fN_{+}$, where $\fN_{\pm} = 
\fn_{\pm}\oplus\fg_{\pm}$, and for the cases when $\fg_{0}$ is a Lie 
algebra (purely even) describe the defining relations.  The relations 
obtained for vectorial Lie superalgebras are not very simple-looking, 
cf.  \cite{LP}.

Notice that, unlike the case of finite dimensional simple Lie algebras, 
the bases, i.e., systems of simple roots, correspond not to maximal 
solvable Lie superalgebras (described in \cite{Shc}) but to what is 
called {\it Borel subalgebras} in \cite{PS}.

{\bf Open problems} are listed in \S 4.

\section* {\protect \S 1.  Generators in some vectorial Lie 
superalgebras and associated loops} 

\ssec{1.1.  Generators of $\fvect (0|n)$} Set $\partial _i = 
\frac{\partial}{\partial x_i}$.  Some of the generators of $\fvect 
(0|n)$ generate its subalgebras as indicated (i.e., the $X^+_{i}$ and 
$Y$ generate $\fN_{+}$; the $X^{\pm}_{i}$ generate $\fsl(1|n)$):
\vskip 0.2 cm

\small
\renewcommand{\arraystretch}{1.4}
\begin{tabular}{|c|c|c|}
\hline
&$\fsl (1|n)$&\\
\hline
$\fN_+$& $x_1\partial _2,  \; \; \;  \dots , x_{n-1}\partial  _n, \; \; \;   
{\bf x_n\sum x_i\partial _i}$ & $x_nx_{n-1}\partial _1$\\
\hline
$\fN_-$&$x_2\partial _1, \; \; \;  \dots , \; \; \;  x_n\partial _{n-1}, \; \;
\; \; \partial  _n$&\\ 
\hline
{\rm notations}& $X^{\pm}_1,\; \; \; \;    \dots , \; \; \;  X^{\pm}_{n-1},    
\; \; \; X^{\pm}_n$&$Y$\\ 
\hline
\end{tabular}
\normalsize
\vskip 0.2 cm

\noindent The generators of $\fsvect (0|n)$ are the same as of $\fvect 
(0|n)$ but without the boldfaced element $X^+_n=x_n\sum x_i\partial 
_i$.  The loop algebras have two more generators: for $\fvect (0|n)$ 
set 
$$
X^-_0=x_1\dots x_n\partial _n\cdot t^{-1}\; \text{ and }\; X^+_0=\partial 
_1\cdot t.
$$
For $\fsvect (0|n)$ set 
$$
X^-_0=x_1\dots x_{n-1}\partial _n\cdot t^{-1}\; \text{ and }\;
X^+_0=\partial _1\cdot t.
$$
It is not clear that this choice of generators (the 
highest weight vector of $\fg_{-1}$ and the lowest weight vector of 
$\fg_1$) which gives nice-looking relations for Lie algebras (and even 
Lie superalgebra with Cartan matrix) is the best when Lie 
superalgebras are very non-symmetric.

\ssec{1.2.  Generators of $\fk (1|n)$} In what follows we will by 
abuse of language write just $f$ instead of either $H_f$, the 
Hamiltonian vector field generated by $f$ or $K_f$, the contact vector 
field generated by $f$; in so doing we must remember that in either 
case ($H_f$ or $K_f$) the degree of the vector field generated by a 
monomial $f$ of degree $k$ is equal to $k-2$.

Some of the generators of $\fk (1|n)$ generate the following 
subalgebras (for $n = 2k>6$ and $n = 2k+1>5$, respectively):

\vskip 0.2 cm

\noindent
\footnotesize \renewcommand{\arraystretch}{1.4}
\begin{tabular}{|c|c|c||c|c|}
\hline
&$\fosp (2k|2)$&&$\fosp (2k+1|2)$&\\
\hline
$\fN _+$ & ${\bf t\eta_1}\; \; \; \;   \xi_1\eta_2 \; \;  \dots \; \; 
\eta_n\xi_{n-1} \; \;  \; \; \xi_n\xi_{n-1} $&$\eta_1\eta_2\eta_3$& ${\bf
t\eta_1}\; \; \; \;   \xi_1\eta_2 \; \;  \dots \; \;  \eta_n\xi_{n-1} \; \; \;
\; \xi_n\theta$&$\eta_1\eta_2\eta_3$\\  
\hline
$\fN _- $& $\xi_1\; \; \; \; \eta_1\xi_2\; \; \;    \dots \; \; \;  
\eta_{n-1}\xi_n\; \; \; \; \eta_n\eta_{n-1}$&&$\xi_1\; \; \; \; \eta_1\xi_2\;
\; \;    \dots \; \; \;   \eta_{n-1}\xi_n\; \; \; \;   \eta_n\theta$&\\ 
\hline
\text{notations} &$X^{\pm}_0\; \; \;  X^{\pm}_1\; \;  \dots \; \;
X^{\pm}_{n-1}\; \; \;  X^{\pm}_n $&$Y$&$X^{\pm}_0\; \; \;  X^{\pm}_1\; \; \dots
\; \; X^{\pm}_{n-1}\; \; \;  X^{\pm}_n $&$Y$\\ 
\hline
\end{tabular}
\normalsize
\vskip 0.2 cm

\noindent The generators of $\fh (0|n)$ and $\fpo (0|n)$ are those 
above without the boldfaced element $X^+_{0}= t\eta_1$.

The loop algebras $\fh (0|n)^{(1)}$ and $\fpo (0|n)^{(1)}$ have two 
more generators:
$$
Y^-_0=\xi_1\dots \xi_n\eta_n\eta_{n-1}\dots\eta_{2}\cdot t^{-1}\;
\text{ and }\; Y^+_0=\eta_1\cdot t.
$$
It is not clear that this choice of generators is the best and it 
is desirable to experiment with other choices.

In small dimensions ($n<7$) relations look differently and are to be 
computed separately.  Besides, the generators look different.  Though 
presentations of some of these algebras were considered, it is 
advisabable to revise it (the results of A.~Nilsson are unpublished 
and those of \cite{FNZ}, as well of \cite{T}, should be presented in a 
more user-friendly form).

\section* {\protect  \S 2. Relations}

\ssec{2.1. Relations for $\fN_-$ of $\fk(1|n)$, $\fpo(0|n)$ and $\fh(0|n)$}
Clearly, $\fN_-$ for $\fvect(0|n)$ and $\fsvect (0|n)$ coincides with $\fn_-$ for
$\fsl(1|n)$ while $\fN_-$ for $\fk (1|n)$  and $\fpo (0|n)$ coincides with
$\fn_-$ for $\fosp(n|2)$, which  are known \cite{LSe}, \cite{GL2}.

\ssec{2.1.1.  Relations for $\fN_-$ of $\fh(0|n)$} The 
Lie algebra $\fh(0|n)$ is generated by the same elements as $\fk(1|n)$ 
and $\fpo(0|n)$ but the relations are different: for $\fh(0|n)$ there 
is an additional relation of weight $(0,\ \dots,\ 0)$ with respect to 
$\fo(n)$ because (for $n>1$)
$$
H_2(\fg_{-1}) = S^2(\fg_{-1}) = 
R(2\pi) \oplus R(0).
$$
The corresponding cycle of weight 0 is 
$$
\{\xi_1, \eta_1\} + \dots + \{\xi_n, \eta_n\} \; \; \; (+ \{\theta, \theta\}\;
{\rm if}\; n\; \text{ is odd }).\eqno{(*)} 
$$
The relation expressed in terms of generators looks awful.  It can be 
beautified as follows.  In the space of relations corresponding to the 
other irreducible component the subspace of relations of weight 0 is 
of dimension $n-1$.  Therefore, $n-1$ summands in $(*)$ vanish; for 
role of survivor select 
the simplest one of them, say, the following one:
$\{\xi_1, \eta_1\} = 0$.

\ssec{2.1.2.  Relations for $\fN_+$ of $\fh (0|n)$, $n>4$}
$H_2(\fg_+)$ is the direct sum of irreducible $\fg_0$-modules with the
following lowest weights with respect to $\fo(n)$, for notations see
Tables in \cite{OV}:

For $n=2l$, $l\geq 5$:
\vskip 0.2 cm

\noindent
\small\renewcommand{\arraystretch}{1.4}
\begin{tabular}{|c|c|c|}
\hline
$N$ & the lowest weight & the corresponding cycle \\
\hline
1 & $-2(\eps _1+\eps _2 +\eps _3)$ & $ 
\eta_1\eta_2\eta_3\wedge \eta_1\eta_2\eta_3 $\\
\hline
2 &$ -2(\eps _1+\eps _2) $ & $\sum 
_i\eta_1\eta_2\eta_i\wedge \eta_1\eta_2\xi_i$\\
\hline 
3 & $-2\eps _1 $ & $ \sum 
_{i,j}\eta_1\eta_i\xi_j\wedge\eta_1\eta_i\xi_j- 2\sum 
_{i<j}\eta_1\eta_i\eta_j\wedge \eta_1\xi_i\xi_j$ \\
\hline
4 &$ 0 $&$  \sum _{i, j}\sum _{k\leq
l}(\eta_i\eta_j\xi_k\wedge\xi_i\xi_j\eta_k+\eta_i\eta_j\eta_k
\wedge\xi_i\xi_j\xi_k)$\\ 
\hline 
5 & $-(2\eps _1+\eps _2 +\eps _3+\eps _4 
+\eps _5)$ & $\eta_1\eta_2\eta_3\wedge \eta_4\eta_5\eta_1 
-\eta_1\eta_2\eta_4\wedge \eta_3\eta_5\eta_1 +
\eta_1\eta_2\eta_5\wedge 
\eta_3\eta_4\eta_1 $ \\
\hline 
6 & $-\eps _1$ &$\eta_1
\sum_i\eta_i\eta_{i+1}\eta_{i+2}\wedge \xi_i\xi_{i+1}\xi_{i+2}$ \\
\hline
\end{tabular}

\normalsize
\vskip 0.2 cm

For small $l$ the relations look differently; the form of relations 1) 
-- 3) is the same as in the general case, the new in form relations 
are (here $\sum _{{\rm cycl}}$ means the cyclic permutation 
of $\eta_1\eta_2\eta_3$):

\noindent
\small\renewcommand{\arraystretch}{1.4}
\begin{tabular}{|c|c|c|c|}
\hline
$l$ & N & the lowest weight & the corresponding cycle \\
\hline
3 &5 & $-2\eps _1$ & $ \eta_1\eta_2\xi_2\wedge
\eta_1\eta_3\xi_3 -\eta_1\eta_2\xi_3\wedge
\eta_1\eta_3\xi_2 -\eta_1\eta_2\eta_3\wedge
\eta_1\xi_2\xi_3 $\\ 
\hline
3 & 6 &$ -2(\eps _3-\eps _1-\eps _2) $ &
$\eta_1\eta_2\xi_3\wedge \eta_1\eta_2\xi_3$\\ 
\hline
\hline 
4 & 5 & $\eps _4-\eps _1+\eps _2 +\eps _3 
$ & $ \sum _{{\rm cycl}}\eta_1\eta_2\xi_4\wedge\eta_3\eta_4\xi_4- \sum 
_{i}\sum _{{\rm cycl}}(\eta_1\eta_2\xi_4\wedge \eta_3\eta_i\xi_i+$ \\
&&&$\eta_1\eta_2\eta_3\wedge\sum _{i}\eta_i\xi_i\xi_4)$\\
\hline 
4& 6 &$-(2\eps _1+\eps _2
+\eps _3)$&$\eta_1\eta_2\eta_3\wedge\eta_1\eta_4\xi_4-
\eta_1\eta_2\eta_4\wedge\eta_1\eta_3\xi_4+
\eta_1\eta_3\eta_4\wedge\eta_1\eta_2\xi_4$\\  
\hline  
\hline 
5 & 6 & $ -(2\eps _1-\eps _2 -\eps _3-\eps 
_4+ \eps _5)$ & $\eta_1\eta_2\eta_3\wedge \eta_4\eta_5\eta_1 
-\eta_1\eta_2\eta_4\wedge \eta_3\eta_5\eta_1 +\eta_1\eta_2\eta_5\wedge 
\eta_3\eta_4\eta_1 $ \\
\hline 
6 & &$-\eps _1$ &$\eta_1
\sum_i(\eta_i\eta_{i+1}\eta_{i+2}\wedge
\xi_i\xi_{i+1}\xi_{i+2})$\\
\hline
\end{tabular}
\normalsize
\vskip 0.2 cm
 
For $n=2l+1$, $l\geq 5$:
\vskip 0.2 cm

\noindent
\small\renewcommand{\arraystretch}{1.4}
\begin{tabular}{|c|c|c|}
\hline
$N$ & the lowest weight & the corresponding cycle \\
\hline
1 & $-2(\eps _1+\eps _2 +\eps _3)$ & $ 
\eta_1\eta_2\eta_3\wedge \eta_1\eta_2\eta_3 $\\
\hline
2 &$ -2(\eps _1+\eps _2) $ & $\sum 
_i\eta_1\eta_2\eta_i\wedge \eta_1\eta_2\xi_i$\\
\hline 
3 & $-2\eps _1 $ &$\sum 
_{i,j}\eta_1\eta_i\xi_j\wedge\eta_1\eta_i\xi_j- 2\sum 
_{i<j}\eta_1\eta_i\eta_j\wedge \eta_1\xi_i\xi_j$ \\
\hline
4 &$ 0 $&$  \sum _{i, j}\sum _{k\leq
l}(\eta_i\eta_j\xi_k\wedge\xi_i\xi_j\eta_k+\eta_i\eta_j\eta_k
\wedge\xi_i\xi_j\xi_k)$\\ 
\hline 
5 & $-(2\eps _1+\eps _2 +\eps _3+\eps
_4+\eps _5)$ &  $\eta_1\eta_2\eta_3\wedge
\eta_4\eta_5\eta_1 -\eta_1\eta_2\eta_4\wedge
\eta_3\eta_5\eta_1 +\eta_1\eta_2\eta_5\wedge
\eta_3\eta_4\eta_1    $ \\
\hline 
6 & $-\eps _1$ &$\eta_1
\sum_i\eta_i\eta_{i+1}\eta_{i+2}\wedge \xi_i\xi_{i+1}\xi_{i+2}$\\
\hline
\end{tabular}
\vskip 0.2 cm

\normalsize

The space $H_1(\fn_+; \fg_1)$ is responsible for the following 
relations (element indicated should vanish for $i>3$):
$$
\{\xi_1\eta_2, \eta_1\eta_2\eta_3\},\;
\{\xi_2\eta_3, \eta_1\eta_2\eta_3\},\;
(\ad~\xi_3\eta_4)^2(\eta_1\eta_2\eta_3),\;
\{\xi_i\eta_{i+1}, \eta_1\eta_2\eta_3\},\;
\{\xi_{n-1}\xi_n, \eta_1\eta_2\eta_3\}.
$$


\ssec{2.1.3.  Relations for $\fN_+$ of $\fk (1|n), n>4$} For the 
$X_{i}^{+}$ the relations are the same as for $\fn_{+}$ of 
$\fosp(n|2)$, cf.  \cite{GL2}.

The relations between the $X_i^+$, $1 \leq i\leq n$ and $Y$ are the 
same as for $\fh (0|n)$.  New relations involving $X_0^+$ and $Y$ are:

\vskip 0.2 cm

\noindent
\small\renewcommand{\arraystretch}{1.4}
\begin{tabular}{|c|c|c|}
\hline
$N$ & the lowest weight & the corresponding cycle \\
\hline
1 & $-4\eps _1 $&$ \sum _i(-q_1^2q_i\wedge q_1^2p_i)
+(n+2)tq_1\wedge q_1^3 $  \\
\hline
2 & $-3\eps _1-\eps _2$ &$ q_1^3 \wedge tq_2
+q_1^2q_2\wedge tq_1 $\\
\hline
\end{tabular}
\vskip 0.2 cm

\normalsize

\ssec{2.1.4.  Relations between $\fN_+$ and $\fN_+$ for $\fh (0|n)$ 
and $\fk (1|n)$, $n>4$} These relations are as for $\fosp(2|2n)$ 
unless they involve $Y$; and the new extra ones are:
$$
[Y, X_0^-] = \eta_2\eta_3\qquad
[Y, X_i^-] = 0\text{ for }\ i>0.
$$

\ssec{ 2.2. Relations for $\fvect(0|n)$ and $\fsvect (0|m)$, $m>2$}

\ssec{2.2.1.  Relations for $\fN_+$ of $\fvect (0|n)$, $n>2$} The 
space $H_1(\fn_+; \fg_1)$ is spanned by $(\ad p_1q_2)^3q_1^3\wedge 
p_1q_2, \quad q_1^3\wedge p_2q_3$, \dots , $q_1^3\wedge 
p_{n-1}q_n$, and $q_1^3\wedge p_n^2$.

$H_2(\fg_+)$ is the direct sum of irreducible $\fg_0$-modules
with the following lowest weights:

\vskip 0.2 cm

\noindent
\small\renewcommand{\arraystretch}{1.4}
\begin{tabular}{|c|c|c|}
\hline
$N$ & the lowest weight & the corresponding cycle \\
\hline
1 &$ 2(\eps _n+\eps _{n-1}-\eps _ 1) $&
$\xi_n\xi_{n-1}\partial_1\wedge
\xi_n\xi_{n-1}\partial_1$\\ 
\hline
2 &$2\eps _n+\eps _{n-1}+\eps _{n-2}-\eps _ 1-\eps _2  $ & $
\xi_n\xi_{n-1}\partial_2\wedge
\xi_n\xi_{n-2}\partial_1-\xi_n\xi_{n-1}\partial_1\wedge
\xi_n\xi_{n-2}\partial_2$ \\  
\hline 
3 &$2\eps _n+\eps _{n-1} -\eps _ 1$ & $ \sum 
_i\xi_n\xi_i\partial_1\wedge \xi_n\xi_{n-1}\partial_i $  \\
\hline
4 &$ \eps _n+\eps _{n-1}+\eps _{n-2}+
\eps _{n-3}-2\eps _1 $&$ 
\xi_n\xi_{n-1}\partial_1\wedge\xi_{n-2}\xi_{n-3}\partial_1- 
\xi_n\xi_{n-2}\partial_1\wedge\xi_{n-1}\xi_{n-3}\partial_1+$\\
&&$\xi_n\xi_{n-3}\partial_1\wedge\xi_{n-1}\xi_{n-2}\partial_1$  \\ 
\hline 
5 & $2\eps _n$ &$ \sum  _{i, j}\xi_n\xi_i\partial_i
\wedge\xi_n\xi_j\partial_j$ \\ 
\hline
6 &$ \eps _n+\eps _{n-1}+\eps _{n-2} - \eps _1 $&$ 
\sum  _i(\xi_n\xi_{n-1}\partial_1\wedge\xi_{n-2}\xi_i\partial_i+ 
\xi_{n-2}\xi_{n-3}\partial_1\wedge\xi_n\xi_i\partial_i+$\\
&&$\xi_{n-2}\xi_n\partial_1\wedge\xi_{n-1}\xi_i\partial_i)$  \\ 
\hline
7 &$ 2\eps _n+\eps _{n-1}-\eps _1 $&$ 
\sum  _i\xi_n\xi_i\partial_i\wedge\xi_n\xi_{n-1}\partial_1$  \\ 
\hline
8 & $2\eps _n$ &$ (\sum   \xi_n\xi_i\partial_i)
\wedge (\sum   \xi_n\xi_j\partial_j)$ \\ 
\hline
\end{tabular}

\normalsize
\vskip 0.2 cm

The corresponding relations for $\fsvect$ are the relations 1) -- 4).

The relations for $\fvect (0|4)$ are
\vskip 0.2 cm

\noindent
\small\renewcommand{\arraystretch}{1.4}
\begin{tabular}{|c|c|c|}
\hline
$N$ & the lowest weight & the corresponding cycle \\
\hline
1 &$ 2(\eps _4+\eps _3-\eps _ 1) $&$\xi_4\xi_3\partial_1\wedge
\xi_4\xi_3\partial_1$\\ 
\hline
2 &$2\eps _4+\eps _3-\eps _ 1$ & $
\sum  \xi_4\xi_i\partial_1\wedge
\xi_4\xi_3\partial_i$ \\  
\hline 
3 &$2\eps _4$ & $ \sum  _{i, j}
\xi_4\xi_i\partial_j\wedge \xi_4\xi_j\partial_i $  \\
\hline
4 &$ \eps _4+\eps _3+\eps _2-\eps _1 $&$ 
\sum  (\xi_4\xi_3\partial_1\wedge\xi_2\xi_i\partial_i+ 
\xi_4\xi_2\partial_1\wedge\xi_3\xi_i\partial_i+
\xi_3\xi_2\partial_1\wedge\xi_4\xi_i\partial_i)$  \\ 
\hline 
5 & $2\eps _4$ &$ (\sum  \xi_4\xi_i\partial_i)
\wedge(\sum  \xi_4\xi_i\partial_i)$ \\ 
\hline
6 &$ 2\eps _4+\eps _3-\eps _1 $&$ 
(\sum  \xi_4\xi_i\partial_i)\wedge\xi_4\xi_3\partial_1$  \\ 
\hline
\end{tabular}
\normalsize
\vskip 0.2 cm

The corresponding relations for $\fsvect (0|4)$ are the relations 1) -- 3).
The relations for $\fvect (0|3)$ are

\vskip 0.2 cm
\noindent
\small\renewcommand{\arraystretch}{1.4}
\begin{tabular}{|c|c|c|}
\hline
$N$ & the lowest weight & the corresponding cycle \\
\hline
1 &$ 2(\eps _3+\eps _2-\eps _ 1) $&
$\xi_3\xi_2\partial_1\wedge
\xi_3\xi_2\partial_1$\\ 
\hline 
2 &$2\eps _3 $ & $ \xi_3\xi_2\partial_1\wedge 
\xi_3\xi_1\partial_2+ \xi_3\xi_1\partial_2\wedge \xi_3\xi_2\partial_1$ 
\\
\hline
3 & $2\eps _3$ &$ (\sum  \xi_3\xi_i\partial_i)
\wedge(\sum  \xi_3\xi_i\partial_i)$ \\ 
\hline
4 &$ 2\eps _3+\eps _2-\eps _1 $&$ 
(\sum  \xi_3\xi_i\partial_i)\wedge\xi_3\xi_2\partial_1$  \\ 
\hline
\end{tabular}
\normalsize
\vskip 0.2 cm

The corresponding relations for $\fsvect(0|3) \simeq\fspe(3)$ are the 
relations 1) -- 2).

\ssec{2.2.2.  Relations between $\fN_+$ and $\fN_-$ for $\fvect 
(0|n)$, $n>3$, and $\fsvect (0|m)$, $m>2$} These relations are as for 
$\fsl(1|n)$ unless they involve $Y$; the extra relations are:
$$
[Y, X_n^-]= (x_{n-1}\partial _1) = [ X_{n-1}^-, [\dots , 
[X_2^-, X_1^-]] \dots ]; \qquad
[Y, X_i^-]  = 0\text{ for }\ i>0.  
$$

\ssec{2.2.3.  Relations for $\fN_\pm$ of $\fvect 
(0|n)^{(1)}$, $n>3$} The new relations that involve $X_0^+$ are (we 
only indicate the terms to be equated to zero):

\vskip 0.2 cm

\underline{For $n=3$, $\fN_+$}: $[X_0^+, X_0^+]$, $[X_0^+, X_2^+]$, $[X_0^+,\,Y]$, 
$[X_1^+,\,[X_0^+,\,X_1^+]]$, 
$[Y,\,[X_0^+,\,X_1^+]]-[X_0^+,\,X_3^+]$, ${{(\ad 
X_1^+)}^2}\,[X_1^+,\,Y]$, $[Y,\,{{(\ad X_1^+)}^2}\,Y]$, 
$[[X_0^+,\,X_1^+],\,[X_0^+,\,X_3^+]]$, 
$[[X_0^+,\,X_3^+],\,[X_2^+,\,X_3^+]]$, 

\noindent $[[X_3^+,\,[X_0^+,\,X_1^+]],\,[Y,\,[X_1^+,\,X_2^+]]]+{\frac{1}{2}}\, 
[{{(\ad X_1^+)}^2}\,Y,\,[X_2^+,\,[X_0^+,\,X_3^+]]]$.

\underline{For $n=3$, $\fN_-$}: $[X_0^-,\,X_1^-]$, ${{(\ad X_0^-)}^2}\,X_2^-$, 
${{(\ad X_0^-)}^2}\,X_3^-$, ${{(\ad X_2^-)}^2}\,X_0^-$, 
$[[X_0^-,\,X_2^-],\,[X_0^-,\,X_3^-]]$, 
$[[X_2^-,\,[X_0^-,\,X_3^-]],\,[X_3^-,\,[X_0^-,\,X_2^-]]]$, 

\noindent $[[[X_0^-,\,X_3^-],\,[X_1^-,\,X_2^-]],\, 
[[X_0^-,\,X_3^-],\,[X_2^-,\,X_3^-]]]+2\, 
[[X_3^-,\,[X_1^-,\,[X_0^-,\,X_2^-]]],\, 
[[X_0^-,\,X_3^-],\,[X_2^-,\,X_3^-]]]$

\underline{For $n=4$, $\fN_+$}: $[X_0^+,\,X_0^+]$, $[X_0^+,\,X_2^+]$, 
$[X_0^+,\,X_3^+]$, $[X_0^+,\,Y]$, ${{(\ad X_1^+)}^2}\,X_0^+$, 
$[Y,\,[X_0^+,\,X_1^+]]$, 

\noindent $[[X_0^+,\,X_1^+],\,[X_0^+,\,X_4^+]]$, 
$[[X_0^+,\,X_1^+],\,[X_2^+,\,Y]]-[X_0^+,\,X_4^+]$, 

\noindent $[[X_0^+,\,X_4^+],\,[X_3^+,\,X_4^+]]$, 
$[[X_2^+,\,Y],\,[X_3^+,\,[X_0^+,\,X_4^+]]]$, 

\noindent $[[X_4^+,\,[X_0^+,\,X_1^+]],\,[Y,\,[X_2^+,\,X_3^+]]]- {\frac{1}{2}}\, 
[[X_3^+,\,[X_0^+,\,X_4^+]],\,[Y,\,[X_1^+,\,X_2^+]] ]$, 

\noindent $[[[X_0^+,\,X_4^+],\,[X_1^+,\,X_2^+]],\, 
[[X_1^+,\,Y],\,[X_2^+,\,X_3^+]]]+\,[[[X_0^+,\,X_1^+],\,[X_3^+,\,X_4^+]],\, 
[[X_1^+,\,X_2^+],\,[X_2^+,\,Y]]]$

\underline{For $n=4$, $\fN_-$}: $[X_0^-,\,X_0^-]$, $[X_0^-,\,X_1^-]$, 
$[X_0^-,\,X_2^-]$, ${{(\ad X_3^-)}^2}\,X_0^-$, 

\noindent $[[X_0^-,\,X_3^-],\,[X_0^-,\,X_4^-]]$, 
$[[X_0^-,\,X_4^-],\,[X_3^-,\,[X_0^-,\,X_4^-]]]$, 

\noindent 
$$
\renewcommand{\arraystretch}{1.4}
\begin{array}{l}
	{}[[X_2^-,\,[X_0^-,\,X_3^-]],\,[X_4^-,\,[X_0^-,\,X_3^-]]], \\
{}[[X_4^-,\,[X_0^-,\,X_3^-]],\, [[X_0^-,\,X_4^-],\,[X_2^-,\,X_3^-]]], \\
{}[[X_0^-,\,X_4^-],\,[X_3^-,\,X_4^-]],\, 
[[X_2^-,\,X_3^-],\,[X_4^-,\,[X_0^-,\,X_3^-]]]], 
\end{array}
$$
\noindent $$
\renewcommand{\arraystretch}{1.4}
\begin{array}{l}
	{}[[[X_1^-,\,X_2^-],\,[[X_0^-,\,X_4^-],\,[X_3^-,\,X_4^-]]],\, 
[[X_4^-,\,[X_0^-,\,X_3^-]],\,[X_4^-,\,[X_2^-,\,X_3^-]]]]+\\
{\frac{1}{2}}\, 
[[[X_1^-,\,X_2^-],\, [[X_0^-,\,X_4^-],\,[X_3^-,\,X_4^-]]],\, 
[[X_3^-,\,[X_0^-,\,X_4^-]],\,[X_4^-,\,[X_2^-,\,X_3^-]]]]
\end{array}
$$

\underline{For $n=5$, $\fN_+$}: $[X_0^+,\,X_0^+]$, $[X_0^+,\,X_2^+]$, 
$[X_0^+,\,X_3^+]$, $[X_0^+,\,X_4^+]$, $[X_0^+,\,Y]$, ${{(\ad 
X_1^+)}^2}\,X_0^+$, $[Y,\,[X_0^+,\,X_1^+]]$, 

\noindent $[[X_0^+,\,X_1^+],\,[X_0^+,\,X_5^+]]$, 
$[[X_0^+,\,X_5^+],\,[X_4^+,\,X_5^+]]$, 

\noindent $[[X_3^+,\,Y],\,[X_2^+,\,[X_0^+,\,X_1^+]]]+[X_0^+,\,X_5^+]$, 
$[[X_3^+,\,Y],\,[X_4^+,\,[X_0^+,\,X_5^+]]]=0$, 

$[[X_5^+,\,[X_0^+,\,X_1^+]],\,[Y,\,[X_3^+,\,X_4^+]]]$, 

\noindent $[[Y,\,[X_3^+,\,X_4^+]],\,[[X_0^+,\,X_5^+],\,[X_1^+,\,X_2^+]]]+ 
\,[[Y,\,[X_1^+,\,X_2^+]],\,[[X_0^+,\,X_5^+],\,[X_3^+,\,X_4^+]]]$ 

\noindent $$
\renewcommand{\arraystretch}{1.4}
\begin{array}{l}
	{}[[[X_3^+,\,Y],\,[X_3^+,\,[X_1^+,\,X_2^+]]],\, 
[[X_4^+,\,X_5^+],\,[X_2^+,\,[X_0^+,\,X_1^+]]]]- \\
{}[[[X_2^+,\,X_3^+],\, [X_5^+,\,[X_0^+,\,X_1^+]]],\, 
[[X_3^+,\,X_4^+],\,[Y,\,[X_1^+,\,X_2^+]]]], 
\end{array}
$$
\noindent 
$$
\renewcommand{\arraystretch}{1.4}
\begin{array}{l}
	{}[[[X_4^+,\,X_5^+],\,[[X_1^+,\,Y],\,[X_2^+,\,X_3^+]]],\, 
[[X_3^+,\,[X_1^+,\,X_2^+]],\,[Y,\,[X_3^+,\,X_4^+]]]]=\\
{\frac{1}{2}}\, 
[[[X_2^+,\,X_3^+],\,[X_4^+,\,X_5^+]],\, 
[[X_4^+,\,X_5^+],\,[X_3^+,\,[X_1^+,\,Y]]]] +\\
{\frac{1}{2}}\, 
[[[X_3^+,\,Y],\,[X_4^+,\,X_5^+]],\, 
[[X_4^+,\,X_5^+],\,[X_3^+,\,[X_1^+,\,X_2^+]]]]\end{array}
$$

\underline{For $n=5$, $\fN_-$}: $[X_0^-,\,X_1^-]$, $[X_0^-,\,X_2^-]$, 
$[X_0^-,\,X_3^-]$, ${{(\ad X_0^-)}^2}\,X_4^-$, ${{(\ad 
X_4^-)}^2}\,X_0^-$, $[[X_0^-,\,X_4^-],\,[X_0^-,\,X_5^-]]$, 

\noindent $[[X_4^-,\,[X_0^-,\,X_5^-]],\,[X_5^-,\,[X_0^-,\,X_4^-]]]$, 
$[[X_5^-,\,[X_0^-,\,X_4^-]],\, [[X_0^-,\,X_5^-],\,[X_3^-,\,X_4^-]]]$, 

\noindent $[[[X_0^-,\,X_5^-],\,[X_3^-,\,X_4^-]],\, 
[[X_0^-,\,X_5^-],\,[X_4^-,\,X_5^-]]]$, 

\noindent $[[[X_2^-,\,X_3^-],\,[X_5^-,\,[X_0^-,\,X_4^-]]],\, 
[[X_4^-,\,X_5^-],\,[X_3^-,\,[X_0^-,\,X_4^-]]]]$, 

\noindent $[[[X_4^-,\,X_5^-],\,[X_3^-,\,[X_0^-,\,X_4^-]]],\, 
[[X_2^-,\,X_3^-],\,[[X_0^-,\,X_5^-],\,[X_4^-,\,X_5^-]]]]$

\ssec{2.2.4.  The periplectic series} Recall that the compatible (with 
parity) $\Zee$-gradings of $\fspe (n)$ are of the form $\fspe 
(n)=\fg_{-1} \oplus \fg_{0} \oplus \fg_{1}$ and there are two such 
cases both with ${\fg}_ 0 =\fsl(n)$: (here $\id$ is the standard $\fsl 
(n)$-module):

a) $\fg_{1}=S^{2}(\id)$, $\quad\fg_{-1}=E^{2}(\id^{*})$;

b) $\fg_{1}=  E^{2}(\id)$, $\quad\fg_{-1}=   S^{2}(\id^{*})$ .

\vskip 0.2 cm

Let $\fn^{\pm}$ be the maximal nilpotent
subalgebras of $\fg_{0}$. Set
$$
\fm^{+} =\fn^{+} \oplus \fg_{1};\quad
\fm^{-} =\fn^{-} \oplus \fg_{-1}
$$
Denote by $X^{+}$ (resp. $X^{-}$) a vector
of lowest (highest) weight in the $\fg_{\rm 0}$-module $\fg_{1}$
(resp. $\fg_{-1}$). The first term $\mathop{\oplus}\limits_{p+q=2}
E_{1}^{p, q}$ of the spectral sequence converging to $H_{ 
2}(\fm^{\pm})$ consists of
$$
E_{1}^{2, 0}=H_{2}(\fn^{\pm}), \quad E_{1}^{1, 1}=
H_{1}(\fn^{\pm}; \fg_{ \pm 1}), \quad E_{1}^{0, 2}=H_{0}(\fn^{\pm};
E^{2}(\fg_{\pm 1})).  
$$
Since we already know $H_{2}(\fn^{\pm})$, we are only interested in 
the other two summands. In case (a), (resp. (b)), $H_{1}(\fn^{\pm}; \fg_{\pm1})$
is the same as for $\fm^{+}$ of $\fsp(2n)$
and for $\fm^{-}$ of $\fo(2n)$
(resp. for $\fm^{-}$ of $\fo(2n)$ and $\fsp(2n)$ of
$\fm^{+}$), we explicitly have:
$$
(\ad X_{n}^{+})^{3}(X^{+})=0, \;\;\;\;
(\ad X_{n}^{-})^{2}(X^{-})=0 \eqno{(a)}
$$
$$
(\ad X_{n}^{+})^{2}(X^{+})=0, \;\; \;\;
(\ad X_{n}^{-})^{3}(X^{-})=0 \eqno{(b)}
$$
with $[X_{i}^{+}, X^{+}]=[X_{i}^{-}, X^{-}]=0$ for $i<n$ in both cases (a) and
(b).

Let $\varphi_{i}$ be the $i$th fundamental weight of $\fg_{0}$, 
$R(\chi)$ the (space of the) irreducible representation with highest 
weight $\chi$.  Now, for the $\fsl (n)$-modules 
$\fg_{1}=R(2\varphi_{1})$ in case (a) and $\fg_{1}= R(\varphi_{2})$ in 
case (b), we have:
$$
\renewcommand{\arraystretch}{1.4}
\begin{array}{ll}
	S^{2}(R(2\varphi_{1}))= R(4\varphi_{1})
\oplus R(2\varphi_{2})&\\
S^{2}(R(\varphi_{2})) = \left\{\begin{matrix}R(2\varphi_{2})
\oplus R(\varphi_{4}) &\text{ if }\; n>3\\
R(2\varphi_{2})&\text{ if }\;  n=3.\end{matrix}\right.\end{array}
$$
Therefore, we have the relations
$$
{}[X^{\pm}, X^{\pm}]=0 \eqno{\text{ for both cases a) and b)}}
$$
and the relations 
$$
\renewcommand{\arraystretch}{1.4}
\begin{array}{lll}
	\text{(a)}
&[X^{+}, [X_{n}^{+}, [X_{n}^{+}, X^{+}]]]=0,
&[X^{-}, [X_{n}^{-}, [X_{n-1}^{-}, X^{-}]]]=0; \\  
\text{(b)}&[X^{+}, [X_{n}^{+}, [X_{n}^{+}, X^{+}]]]=0, &
{}[X^{-}, [X_{n}^{-}, [X_{n-1}^{-}, X^{-}]]]=0.
\end{array} \eqno{(PeR_{\pm})}
$$
of which the first in case (a) and the second in case (b) are only 
defined if $n>3$.
 
\ssec{2.3.1 Exceptional loop algebras: $\fd(\eps)^{(3)}$} Let $\eps$ be 
a primitive cubic root of 1 and $\fd(\eps)$ the deform of the Lie 
superalgebra $\fosp(4|2)$ corresponding to the value of parameter 
equal to $\eps$, i.e., $\fd(\eps)=\fg(A)$ for any of the following Cartan 
matrices (cf. \cite{GL2}):
$$
\begin{pmatrix}0&-1&\eps^2\\
1&0&\eps\\
\eps^2&-\eps&0\end{pmatrix},\text{ or }\begin{pmatrix}2&-1&0\\
\eps&0&\eps^2\\
0&-1&2\end{pmatrix},\text{ or }\begin{pmatrix}2&-1&0\\
1&0&\eps\\
0&-1&2\end{pmatrix},\text{ or }\begin{pmatrix}2&-1&0\\
1&0&\eps^2\\
0&-1&2\end{pmatrix}. 
$$
The algebra $\fd(\eps)$ has an outer automorphism of order 3; select 
the generators of the maximal nilpotent subalgebras of 
$\fd(\eps)^{(3)}$ as follows.  Let $X^{\pm}_1$, $X^{\pm}_2$, 
$X^{\pm}_3$ be the Chevalley generators of $\fd(\eps)$.  Set
$$
\begin{matrix}
Y^{+}_1=\eps X^{+}_1+\eps^2X^{+}_2+X^{+}_3, &Y^{+}_2=[X^{-}_3, 
[X^{-}_1, X^{-}_2]]; \\
Y^{-}_1=\eps X^{-}_1+\eps^2X^{-}_2+X^{-}_3, &Y^{+}_2=[X^{+}_1,
X^{+}_2]-[X^{+}_1, X^{+}_3]+[X^{+}_2, X^{+}_3]. \end{matrix}
$$

The relations between these generators are:
$$
\begin{matrix}
[Y^{+}_2, Y^{+}_2]=0, \\
{}[[Y^{+}_1, Y^{+}_2], [Y^{+}_2, [Y^{+}_1, Y^{+}_1]]]=0,\\
(\ad [Y^{+}_1, Y^{+}_1])^3 ([Y^{+}_1, Y^{+}_2])=0,\\
{}[[Y^{+}_1, Y^{+}_1], [[Y^{+}_2, [Y^{+}_1, Y^{+}_1]], [[Y^{+}_2, 
[Y^{+}_1, Y^{+}_1]], [[Y^{+}_1, Y^{+}_1], [Y^{+}_2, [Y^{+}_1, 
Y^{+}_1]]]]]]=64 Y^{+}_1,\\
{}[[Y^{+}_1, Y^{+}_2], [[Y^{+}_2, [Y^{+}_1, Y^{+}_1]], [[Y^{+}_2, [Y^{+}_1,
Y^{+}_1]], [[Y^{+}_1, Y^{+}_1], [Y^{+}_2, [Y^{+}_1, Y^{+}_1]]]]]]=-96
Y^{+}_2. 
\end{matrix}
$$
$$
\begin{matrix}
(\ad Y^{-}_2)^3 Y^{-}_1=0, \\
{}[[Y^{-}_1, Y^{-}_2], [Y^{-}_1, Y^{-}_2]]-2[Y^{-}_2, [Y^{-}_2, [Y^{-}_1,
[Y^{-}_1]]=0,\\ 
(\ad [Y^{-}_1, Y^{-}_1])^2 [Y^{-}_1, Y^{-}_2]=0,\\
{}[[Y^{-}_1, Y^{-}_1], [[Y^{-}_2, [Y^{-}_1, Y^{-}_1]], [[Y^{-}_1,
Y^{-}_2], [Y^{-}_2, [Y^{-}_1, Y^{-}_1]]]]]=-64 Y^{-}_1,\\ 
{}[[Y^{-}_1, Y^{-}_2], [[Y^{-}_2, [Y^{-}_1, Y^{-}_1]], [[Y^{-}_1,
Y^{-}_2], [Y^{-}_2, [Y^{-}_1, Y^{-}_1]]]]]=-64Y^{-}_2. 
\end{matrix}
$$

\ssec{2.3.2.  Stringy superalgebras} For the exceptional stringy 
superalgebra $\fk\fas^{L}$ introduced in \cite{GLS1} (and in \cite{CK}) the 
relations are computed in \cite{GLS1} and, in another form, in \cite{CK}.  
Observe that the relations of $\fk\fas$ are {\it not} simply the 
relations for $\fk\fas^{L}$ that do not involve the extra generators.

Among the stringy superalgebras $\fk^L(1|6)$ is one of the most 
interesting: it possesses a nondegenerate invariant symmetric bilinear 
form and, therefore, can be $q$-quantized, cf. \cite{LSp}.

The basis of $\fg_+$ for the standard $\Zee$-grading of $\fk^L(1|6)$:
$$
\begin{matrix}
X^+_1=\xi_1\eta_2, &X^+_2=\xi_2\eta_3, &X^+_3=\xi_2\xi_3,
&X^+_0=t\eta_1,\\ 
\tilde X^+_1=\frac{1}{t}\xi_1\xi_3\eta_2\eta_3, &
\tilde X^+_2=\frac{1}{t}\xi_1\xi_2\eta_1\eta_3, &
\tilde X^+_3=\frac{1}{t}\xi_1\xi_2\xi_3\eta_1, &
\tilde X^+_0=\eta_1\eta_2\eta_3
\end{matrix}
$$
The generators $X^{\pm}_i$ for $i\neq 0$, clearly, generate $\fo (6)$ 
while all of them generate $\fosp(6|2)$.  One expects the same 
relations between them, but the other generators interfere and the 
final result is as follows (we skip the superscript):
$$
[X_2, X_3]=0,\; \; [X_2, X_0]=0, \; \; [X_3, X_0]=0, \; \; [X_0, X_0]=0;
$$
$$
(\ad X_1)^2X_2=0; \;  \;(\ad X_1)^2X_3=0, \; \;(\ad X_1)^2X_0=0, \; \;
(\ad X_2)^2X_1=0, \; \; (\ad X_3)^2X_1=0;
$$
$$
\begin{matrix}
[X_1, \tilde X_0]=0,\; \; [X_1, \tilde X_1]=0, \; \; [X_2, \tilde X_0]=0,
\; \; [X_2, \tilde X_2]=0,\\
{}[X_2, \tilde X_3]=0,\; \; [X_3, \tilde X_2]=0, \; \; [X_3, \tilde X_3]=0,
\; \; [X_0, \tilde X_0]=0\\
{}[\tilde X_0, \tilde X_0]=0, \; \; [\tilde X_0, \tilde X_1]=0,\; \; [\tilde X_0, \tilde
X_2]=0, \\
{}[\tilde X_1, \tilde X_2]=0, \; \; [\tilde X_1, \tilde X_3]=0, \; \; [\tilde X_2,
\tilde X_3]=0;
\end{matrix}
$$
$$
[X_2, \tilde X_1]+[X_1, \tilde X_2]=0, \; \; [X_3, \tilde X_1]+[X_1, \tilde
X_3]=0;
$$
$$
(\ad X_3)^2\tilde X_0-2[X_0, \tilde X_3]=0;
$$
$$
\begin{matrix}
[\tilde X_0, [X_0, \tilde X_3]]=0,\; \; [\tilde X_1, [X_1, X_2]]=0, \; \;
[\tilde X_1, [X_1, X_3]]=0, \; \; [\tilde X_1, [X_1, X_0]]=0,\\
{}[\tilde X_2, [X_1, X_2]]=0,\; \; [\tilde X_3, [X_1, X_3]]=0, \\
{}[\tilde X_3, [X_3, \tilde X_0]]=0, \; \;[\tilde X_3, [X_0, \tilde X_2]]=0;
\end{matrix}
$$
$$
(\ad \tilde X_1)^2 X_0=0, \; \; (\ad \tilde X_2)^2 X_0=0, \; \; (\ad \tilde
X_3)^2 X_0=0, \; \;  (\ad \tilde X_3)^2 \tilde X_0=0;
$$
$$
[[X_1,  X_2], [X_3, \tilde X_0]]+2[\tilde X_2, [X_0, \tilde X_1]]
+[X_0, [X_1, \tilde X_2]]+[\tilde X_2, [X_1, X_0]]=0;
$$
$$
[[X_1,  X_2], [X_0, \tilde X_2]]=0,\; [[X_1,  X_3], [X_0, \tilde X_3]]=0,\;  
[[X_1,  \tilde X_2], [X_0, \tilde X_2]]=0,\; [[X_1,  \tilde X_3], [X_0, \tilde
X_3]]=0,\; 
$$
$$
[[X_1, \tilde X_2], [X_3, \tilde X_0]]-[\tilde X_2, [X_0, \tilde X_1]]=0,\; 
[[X_0, \tilde X_3], [\tilde X_2, [X_0, \tilde X_1]]=0,\; 
$$
$$
[[X_0, \tilde X_3], [\tilde X_2, [X_1, X_0]]+
[[X_0, \tilde X_3], [\tilde X_0, [X_1, \tilde X_2]]=0,\; 
$$
$$
2[[X_1, \tilde X_2], [X_0, \tilde X_3]]+
2[\tilde X_3, [X_0, [X_1, X_2]]]+
[[X_1, X_3], [X_0, \tilde X_2]]-
[[X_1, X_2], [X_0, \tilde X_3]]=0; 
$$
$$
2[[X_1, \tilde X_3], [X_0, \tilde X_2]]-
2[\tilde X_3, [X_0, [X_1, X_2]]]+
[[X_1, X_2], [X_0, \tilde X_3]]-
[[X_1, X_3], [X_0, \tilde X_2]]=0. 
$$



\section* {\protect \S 4.  Presentation of vectorial Lie 
superalgebras: an overview of open problems}

\ssec{4.1. Vectorial Lie superalgebras with polynomial coefficients and 
loops with values in them} 

\begin{Lemma} {\em (\cite{ShV})} There are embeddings of Lie superalgebras such that 
their negative parts as $\Zee$-graded Lie superalgebras coincide:
$$
\renewcommand{\arraystretch}{1.4}
\begin{array}{l}
\fsl (n+1|m) \subset \fvect (n|m),\quad
\fosp (m|2n+2) \subset \fk (2n+1|m),\\
\fpe (n+1) \subset \fm (n), \quad\fspe (n+1) \subset \fsm (n).
\end{array}
$$
\end{Lemma}
This is already something: it remains to establish presentation of
${\fN}^{+}$ only.  The rest of the known information is gathered in
the following statement, cf.  \cite{HP}.  The cases when the Lie
superalgebra is not simple are marked by a $+$ sign; we disregard
them in this paper.

\begin{Theorem} The following table lists the restrictions for the 
relations of the subalgebra ${\fN}^{+}$ of a simple vectorial Lie 
superalgebra to be ``pure'' {\em (blank space)} or ``dirty'' {\em (marked by a 
dot)}.
\vskip 0.2 cm

\small
$$
\begin{matrix}
\fvect(n|m)&\fsvect(n|m)\\
\hbox{
\tiny
\begin{tabular}{|c|c|c|c|cccc}
$m$&$\;\; \;\;$
&$\;\;\;\;$&$\;\;\;$&$\;\;\;$&$\;\;$
\\
$\;\;$&$\;\;$&$\;\;$&$\;\;$&$\;\;$&$\;\;$
\\
$\;\;$&$\;+ \;$&$\;\;$&$\;\;$&$\;\;$&$\;\;$
&$\;\;$&$\;\;$\\ 
\hline  
3&$\;+ \;$&$\;\;$&$\;\;$&$\;\;$&$\;\;$
&$\;\;$&$\;\;$ \\
\hline 
2&$\;+ \;$&$\;\;$&$\;\;$&$\;\;$&$\;\;$&$\;\;$
&$\;\;$ \\ 
\hline 
1 &$\;+ \;$&$\;\mbox{$\bullet$}\;$&$\;\;$&$\;\;$&$\;\;$
&$\;\;$\\ 
\hline 
0 &$\;+ \;$&$ \mbox{$\bullet$}$
&$ \mbox{$\bullet$}$ &$\;\;$&$\;\;$
&$\;\;$\\ 
\hline 
$\;\;$  & 0 & 1 &2 &3 &$\;\;$&$\; \ldots\;$& $n$ \\ 
\hline 
\end{tabular}}&
\hbox{\tiny
\begin{tabular}{|c|c|c|c|cccc}
$m$&$ \;\;$&$\;\;$&$\;\;$&$\;\;$&$\;\;$\\
$\;\;$&$\;+ \;$&$\;+\;$&$\;\;$&$\;\;$&$\;\;$\\
$\;\;$&$\;+ \;$&$\;+\;$&$\;\;$&$\;\;$&$\;\;$&$\;\;$&$\;\;$\\ 
\hline  
3&$\;+ \;$&$\;+\;$&$\;\;$&$\;\;$&$\;\;$
&$\;\;$&$\;\;$\\
\hline 
2&$\;+ \;$&$\;+ \;$&$\;\;$&$\;\;$&$\;\;$&$\;\;$
&$\;\;$ \\ 
\hline 
1 & +&$\;+\;$&$\;\;$&$\;\;$&$\;\;$&$\;\;$\\ 
\hline 
0& + &$ \mbox{$+$}$&$ \mbox{$\bullet$}$ &$\;\;$&$\;\;$&$\;\;$&$\;\; $\\ 
\hline 
$\;\;$  &0 & 1 & 2 &3 &$\;\;$&$\;  \ldots\;$& $n$ \\ 
\hline
\end{tabular}}\\
&\\&\\
\fh(2n|m)&\fk(2n+1|m)\\
\hbox{
\tiny
\begin{tabular}{|c|c|c|c|cccc}
$m$&$ \;\;$&$\;\;$&$\;\;$&$\;\;$&$\;\;$ & \\
$\;\;$&$\;+\;$&$\;\;$&$\;\;$&$\;\;$&$\;\;$&\\
$\;\;$&$\;+\;$&$\;\;$&$\;\;$&$\;\;$&$\;\;$&$\;\;$\\ 
\hline  3&$\;+\;$&$\;\;$&$\;\;$&$\;\;$&$\;\;$&$\;\;$&$\;\;$ \\ 
\hline 
2&$\;+\;$&$\;\;$&$\;\;$&$\;\;$&$\;\;$&$\;\;$&$\;\;$ \\ 
\hline 
1 & +&$\;\mbox{$\bullet$}\;$&$\;\;$&$\;\;$&$\;\;$
&$\;\;$\\ 
\hline 
0 & + &$ \mbox{$\bullet$}$&$ \mbox{$\bullet$}$ &$\;\;$&$\;\;$&$\;\;$
&$\;\; $\\ 
\hline 
$\;\;$  & 0 & 1 &2 &3 &$\;\;$&$\;  \ldots\;$& $n$ \\ 
\hline
\end{tabular}}&
\hbox{
\tiny
\begin{tabular}{|c|c|c|c|cccc}
$m$&$ \;\;$&$\;\;$&$\;\;$&$\;\;$&$\;\;$\\
$\;\;$&$\;\;$&$\;\;$&$\;\;$&$\;\;$&$\;\;$\\
$\;\;$&$\;\;$&$\;\;$&$\;\;$&$\;\;$&$\;\;$&$\;\;$&$\;\;$\\ 
\hline  
3&$\;\;$&$\;\;$&$\;\;$&$\;\;$&$\;\;$&$\;\;$&$\;\;$ \\
\hline 
2&$\;\bullet\;$&$\;\;$&$\;\;$&$\;\;$&$\;\;$&$\;\;$
&$\;\;$ \\
\hline 
1 & $ \bullet$&$\;\;$&$\;\;$&$\;\;$&$\;\;$&$\;\;$\\ 
\hline 
0& $\bullet$ &$ \mbox{$\bullet$}$&$\;\;$ &$\;\;$&$\;\;$&$\;\;$&$\;\; $\\ 
\hline
$\;\;$  & 0 & 1 &2 &3 &$\;\;$&$\;  \ldots\;$& $n$ \\ 
\hline 
\end{tabular}}\\
&\\&\\
\fsvect^{\circ}(n)&\fsh(n)\\
\hbox{\tiny
\begin{tabular}{|c|c|c|c|ccc} 
\hline  
+&+ &$\; \bullet\; $&  & &  &  \\ 
\hline 
0 & 1 & 2 & 3 &4 &$\ldots$  & $n$ \\
\hline 
\end{tabular}}&
\hbox{
\tiny
\begin{tabular}{|c|c|c|c|ccc} 
\hline
+ &+ &+&+ &  &  &  \\ 
\hline 
0 & 1 & 2 & 3 &4&$\ldots$  & $n$\\ 
\hline 
\end{tabular}}\\
&\\&\\
\fm(n)&\fb_{\lambda}(n), \; \fle(n), \; 
\fsle(n)\\
\hbox{
\tiny
\begin{tabular}{|c|c|c|c|ccc} 
\hline
$\;\bullet\;$ & &  & &  &  &  \\ 
\hline 
0 & 1
& 2 & 3 &4 &$\ldots$  & $n$\\ 
\hline 
\end{tabular}}&
\hbox{
\tiny
\begin{tabular}{|c|c|c|c|ccc} 
\hline
+ &$\bullet$  &  & $\;\;$ &&  &  \\ 
\hline 0 & 1 & 2 & 3 &4 &$\ldots$  & $n$\\
\hline \end{tabular}}
\\
&\\&\\
\fvect(0|n)^{(1)}, \; \fsh(2n)^{(2)}&
\fsvect(0|n)^{(1)}\\
\hbox{
\tiny
\begin{tabular}{|c|c|c|c|ccc}
\hline
+ & $\;\;$ &$ \;\; $& $\;\;$ &&  &  \\ 
\hline 
0 & 1 & 2 & 3 &4 &$\ldots$  & $n$\\
\hline 
\end{tabular}}&
\hbox{\tiny
\begin{tabular}{|c|c|c|c|ccc} 
\hline
+ & + &  & &&  &  \\ 
\hline 
0 & 1 & 2 & 3 &4 &$\ldots$  & $n$\\
\hline 
\end{tabular}}\\
&\\&\\
\fsh(n)^{(1)}&\\
\hbox{\tiny
\begin{tabular}{|c|c|c|c|ccc}
\hline 
+&+&+& & & & \\ 
\hline 
0&1&2&3 &4&$\ldots$  & $n$\\ 
\hline 
\end{tabular}}&\end{matrix}
$$
\end{Theorem} \normalsize

\ssec{4.2.  Twisted loops} There remains to compute the relations in
the following twisted loop superalgebras which have no Cartan matrix:
see \cite{LSS}: $\fpsl(n|n)^{(2)}_{\Pi}$; $\fpsl(n|n)^{(2)}_{\Pi\circ
(-st)}$; $\fpsq(n)^{(4)}_{\Pi}$.

\ssec{4.3.  Stringy and vectorial superalgebras: other bases} We have
only computed the relations in several cases.  There remains to
compute the relations in the cases marked by a thick dot (and a cross)
in Table 4.1 and for the exceptional algebras \cite{Sh5} (for the
consistent gradings thereof this is done in \cite{GLS2}; for the other
W-gradings it is desirable, at least, for completeness).  It is
desirable to compute them for all simple stringy superalgebras or, at
least, for the distinguished algebras that admit nontrivial central
extensions most often used in applications, cf.  \cite{GLS1}, and for
other bases of vectorial algebras and superalgebras.

\end{document}